\documentclass{amsart}
\title{\textbf{On the Beem-Nair Conjecture}}
\author{Shun Furihata}
\date {}
\usepackage{amssymb}
\usepackage{amsmath}
\usepackage{stmaryrd}
\usepackage{graphics}
\usepackage{enumerate}
\usepackage{latexsym}
\usepackage{amsthm}
\usepackage{comment}
\usepackage{pb-diagram}
\usepackage{color}
\usepackage{indentfirst}
\usepackage{titlesec}
\titleformat{\section}{\Large\bfseries}{\thesection}{5pt}{}
\titleformat{\subsection}{\large\bfseries}{\thesubsection}{5pt}{}
\everymath{\displaystyle}
\newcommand{\Ocal}{\mathcal{O}}
\newcommand{\Dch}{\mathcal{D}_{G}^{ch}}
\newcommand{\I}{\mathbf{I}_G}
\newcommand{\g}{\mathfrak{g}}
\newcommand{\n}{\mathfrak{n}}
\newcommand{\tf}{\mathfrak{t}}
\renewcommand{\a}{\widetilde{a}}
\renewcommand{\b}{\widetilde{b}}
\renewcommand{\c}{\widetilde{c}}
\renewcommand{\d}{\widetilde{d}}
\newcommand{\f}{\widetilde{f}}
\newcommand{\gh}{\widehat{\mathfrak{g}}}
\newcommand{\npm}{\mathfrak{n}[t^{\pm1}]}
\newcommand{\Der}{\mathrm{Der}_{\mathbb{C}}}
\newcommand{\DerOG}{\mathrm{Der}_{\mathbb{C}}(\Ocal(G))}
\newcommand{\vac}{\lvert 0 \rangle}
\newcommand{\gr}{\mathrm{gr}}
\newcommand{\J}{J_{\infty}}
\newcommand{\relmiddle}[1]{\mathrel{}\middle#1\mathrel{}}
\newcommand{\bx}{\beta_x}
\newcommand{\by}{\beta_y}
\newcommand{\gx}{\gamma_x}
\newcommand{\gy}{\gamma_y}
\newcommand{\lp}{\lambda_+}
\newcommand{\lm}{\lambda_-}
\newcommand{\bgx}{\langle \bx, \gx \rangle}
\newcommand{\bgy}{\langle \by, \gy \rangle}
\newcommand{\el}{\langle \eta, \lambda_{\pm} \rangle}

\newtheorem{thm}{Theorem}
\newtheorem{lem}[thm]{Lemma}
\newtheorem{prop}[thm]{Proposition}

\theoremstyle{definition}
\newtheorem{defi}[thm]{Definition}
\theoremstyle{definition}
\newtheorem{rmk}[thm]{Remark}
\theoremstyle{remark}
\newtheorem*{prf}{Proof}

\begin{document}
\maketitle

\begin{abstract}
For a simple linear algebraic group $G$, the chiral universal centralizer $\mathbf{I}_{G,k}$ is a vertex operator algebra, which is 
the chiralization of the universal centralizer $\mathfrak{Z}_G$.
The variety $\mathfrak{Z}_G$ is
identified with the spectrum of the equivariant Borel-Moore homology of the affine Grassmannian of the Langlands dual group of $G$.
Beem and Nair conjectured that an open symplectic immersion from $\mathrm{KT}_G$, the Kostant-Toda lattice associated to a simple group $G$, to $\mathfrak{Z}_G$ gives rises to
a free field realization of the chiral universal centralizer at the critical level. 
In this paper, we construct a free field realization of $\mathbf{I}_{G,k}$
at any level, which coincides with the one conjectured by Beem and Nair at the critical level.
We give an explicit description of this construction in $SL_2(\mathbb{C})$-case.
\end{abstract}

\setcounter{section}{-1}
\section*{Introduction}
Vertex operator algebra (VOA) was defined in 1980's (\cite{B}) so as to describe the two-dimensional conformal field theory (CFT) in physics in the algebraic framework.
The study of VOAs may provide new discoveries in physics, and conversely, the research in physics can stimulate the development of new theories of VOAs.
One of the physical aspects in pure mathematics via VOAs is the Higgs Branch Conjecture(\cite{BR, A3}), 
which implies that the Higgs branches (a certain geometrical invariant in physics) are the associated varieties of VOAs (algebraic objects in mathematics) corresponding to four-dimensional $\mathcal{N}=2$ super CFT
via the 4$d$/2$d$ correspondence (\cite{BLL+}).
It is meaningful to find free field realization of VOAs, because the free field realization of the VOA arising from the four-dimensional super CFT 
mirrors the effective field theory description of the super CFT on the Higgs branch of its moduli space of vacua (\cite{BMR}).

The universal centralizer $\mathfrak{Z}_G$ is a family of pairs $(g, s)$, where $s$ is an element in the Slodowy slice and $g \in G$ is in the stabilizer group of $s$.
This is a subvariety of the equivariant Slodowy slice with a symplectic structure(\cite{K3}),
and identified with the space associated to the sphere in the Moore-Tachikawa two-dimensional topological quantum field theory that describes the Higgs branches for the class $\mathcal{S}$ theory(\cite{MT}).
We have an open immersion from a complexification of the phase space of the Kostant-Toda lattice to $\mathfrak{Z}_G$
(\cite{C, BN}).
Also, $\mathfrak{Z}_G$ is identified with the spectrum of the equivariant Borel-Moore homology of the affine Grassmannian (\cite{BFM}) of the Langlands dual group of $G$.
According to \cite{BPRvR}, the chiral counterpart of the universal centralizer $\mathfrak{Z}_G$ is the chiral universal centralizer $\I$ defined in \cite{A1}. 
That is, $\mathfrak{Z}_G$ is the associated variety of $\I$.
For $G = SL_2(\mathbb{C})$ and a generic level, the universal centralizer coincides with the modified regular representation of the Virasoro algebra defined in \cite{FS}.
The representation of the universal centralizer is studied in \cite{MY} for $G = SL_2(\mathbb{C})$ and at level $-1$, but for other general case, it is yet to be studied.
 
In \cite{BN}, the authors suggested how to construct a free field realization of a chiral universal centralizer at the critical level for simple Lie groups $G$.
They conjectured that the immersion $\mathrm{KT}_G \rightarrow \mathfrak{Z}_G$ gives rise to a vertex operator algebra homomorphism between the chiralization of $\mathrm{KT}_G$ and the chiralization of $\mathfrak{Z}_G$.
In this paper, %
we will construct a free field realization of the chiral universal centralizers at any level, using an idea of \cite{FS}.
The universal centralizers are defined by performing the Kostant reduction twice to the cotangent space $T^*G$, 
and the chiral universal centralizers are defined analogously: by performing the quantized Drinfeld-Sokolov reduction of the chiral differential operator algebra $\Dch(G)$ on $G$ twice (\cite{MSV, BD, AG}).
We take the big cell $U$ of the Bruhat decomposition of $G$, and get the free field realization by performing the quantized Drinfeld-Sokolov reduction to restriction map $\Dch(G) \rightarrow \Dch(U)$ twice (Theorem \ref{theorem}). 

\subsubsection*{Acknowledgements.}

Some of the results in this paper appear in the author's master thesis.
The author is grateful to his advisor Professor Tomoyuki Arakawa and Ana Kontrec for valuable discussions and useful comments.

\section{Notations}

In this paper, we let $G$ be a simple affine algebraic Lie group unless otherwise stated. 
Fix a Borel subgroup $B \subset G$ and a splitting $B = N \rtimes T$, where $N \subset G$ is a maximal unipotent subgroup and $T \subset G$ is a maximal torus.
Set $\g := \mathrm{Lie}(G), \, \n := \mathrm{Lie}(N), \, \tf := \mathrm{Lie}(T)$ and let $\Delta$ denote the set of simple roots of $\g$.

Choose a root basis $\{ e^{\alpha} \}_{\alpha \in \Delta}$ of $\g$ with $\alpha \in \g_{\alpha}$, where $\g_{\alpha}$ is the root space.
We fix a principal nilpotent element $\textstyle{f_0 := \sum_{\alpha \in \Delta}e^{-\alpha} \in \n}$ and let $\chi_0 \in \n^* \subset \g^*$ be the Killing dual element of $f_0$.
Note that $\chi_0$ is a character of $\n$.

We regard complex numbers as invariant bilinear forms on $\g$ by considering the bilinear form
$\displaystyle{
\frac{k}{2h^{\vee}}\times (\cdot , \cdot)
}$
for $k \in \mathbb{C}$, where $h^{\vee}$ is the dual Coxeter number of $\g$ and $(\cdot, \cdot)$ is the Killing form on $\g$.


\section{Chiral Universal Centralizer}
\subsection{Kostant reduction}
Let $X$ be an affine Poisson algebraic variety with a Hamiltonian $G$-action with the moment map $\mu_X$.
Then we have a moment map $\mu_N$ for $N$-action;
$$\mu_N : X \rightarrow \g^* \xrightarrow{\text{restriction}} \n^*.$$ 

By the Jacobson-Morozov theorem, $f_0$ can be completed to an $\mathfrak{sl}_2$-triple $(e_0, h_0, f_0)$. 
Let $\g^{e_0}$ be an $\mathrm{ad}_{e_0}$-invariant subspace of $\g$ and $S_{f_0} := f_0 + \g^{e_0} \subset \g$ be the Kostant-Slodowy slice.
\begin{thm}[\cite{K1}]
Assume that $\chi_0$ is a regular value of $\mu_N$ and $N$-action on $\mu_N^{-1}(\chi_0)$ is free.
Then
$$\mu_N^{-1}(\chi_0) /N \cong X \times_{\g^*} S_{f_0}$$
as varieties.
\end{thm}

We call $\mu_N^{-1}(\chi_0) / N$ \textbf{the Kostant reduction} of $X$ with respect to $\mu_X$ and denote $X /\!/_{\chi_0}N.$


\subsection{(Classical) universal centralizer}
We have two commuting actions $\rho_L$ and $\rho_R$ of $G$ on $T^*G \cong G \times \g^*$:
\begin{align*}
&\rho_L(g)(h,\phi) = (hg^{-1}, \mathrm{Ad}^*_g(\phi)), \\
&\rho_R(g)(h,\phi) = (gh,\phi) 
\end{align*}
for $g,h \in G$ and $\phi \in \g^*$, where $\mathrm{Ad}^*$ means the coadjoint action of $G$ on $\g^*$.
In fact, $\rho_L$ and $\rho_R$ are Hamiltonian actions, which correspond to the following moment maps:
\begin{align*}
&\mu_L(h,\phi)=\phi ,\\
&\mu_R(h,\phi)=\mathrm{Ad}^*_h(\phi).
\end{align*}

We perform the Kostant reduction to $T^*G$ with respect to $\mu_L$ and get
\begin{align*}
T^*G /\!/_{\chi_0} N &\cong T^*G \times_{\g^*} S_{f_0} \\
&\cong G \times S_{f_0}.
\end{align*}
Since $\rho_L$ and $\rho_R$ commute with each other, $G \times S_{f_0}$ inherits a Hamiltonian $G$-action $\rho'_R$, corresponding to a new moment map $\mu'_R.$
\begin{align*}
&\rho'_R(g)(h,s) = (gh,s), \\
&\mu'_R(h, s)=\mathrm{Ad}^*_h(\chi_s),
\end{align*}
where $\chi_s$ is the Killing dual of $s \in S_{f_0}.$
Therefore, we can take the Kostant reduction of $G \times S_{f_0}$ with respect to $\mu'_R$ and we get 
$\mathfrak{Z}_G := \{(g, s) \in G \times S_{f_0} \mid \mathrm{Ad}_g(s)=s \}$, called
\textbf{the (classical) universal centralizer} of $G$. 
That is, 
\begin{align*}
(G \times S_{f_0})/\!/_{\chi_0} N &\cong (G \times S_{f_0}) \times_{\g^*} S_{f_0} \\
&\cong \{(g, s) \in G \times S_{f_0} \mid \mathrm{Ad}_g(s)=s \} \\
&= \mathfrak{Z}_G.
\end{align*}

We will follow \cite{BN} so as to obtain a ``classical free field realization''. 
Define $\textstyle{\mathrm{KT}_{\g} := \chi_0 + \tf^* + \sum_{\alpha \in \Delta} \g^*_{\alpha} \setminus \{ 0 \}} \subset \g^*.$
Then $\mathrm{KT}_{\g}$ is a symplectic variety (\cite{K2}) and  it is not difficult to see that $\mathrm{KT}_{\g} \cong T^*(\mathbb{C}^{\times})^{\mathrm{rk} \g}$ as a symplectic variety.
We denote the projection map $\mathrm{pr}  : \mathrm{KT}_{\g} \rightarrow (\mathbb{C}^{\times})^{\mathrm{rk} \g}$.
Let $Z_G : T \rightarrow (\mathbb{C}^{\times})^{\mathrm{rk} \g}$ be the map 
$$t \mapsto (\alpha(t))_{\alpha \in \Delta},$$
where $\alpha(t)$ is defined by $\mathrm{Ad}_t(e^{\alpha})=\alpha(t)e^{\alpha}$ for a simple root $\alpha \in \Delta$.

We define the symplectic variety $\mathrm{KT}_G$ to be the pullback of $\mathrm{pr} : \mathrm{KT}_{\g} \rightarrow (\mathbb{C}^{\times})^{\mathrm{rk} \g}$ and $Z_G : T \rightarrow (\mathbb{C}^{\times})^{\mathrm{rk} \g}.$

$$
\begin{diagram}
\node{\mathrm{KT}_G} \arrow{e,t}{\pi_G} \arrow{s,t}{\pi_T} \node{\mathrm{KT}_{\g}} \arrow{s,r}{\mathrm{pr}} \\
\node{T} \arrow{e,t}{Z_G} \node{(\mathbb{C}^{\times})^{\mathrm{rk} \g}}
\end{diagram}
$$
Since $Z_G$ is \'{e}tale, the same holds for $\pi_G$.
Hence, we can pullback the symplectic form on $\mathrm{KT}_{\g}$ to a symplectic form on $\mathrm{KT}_G.$

We have the following proposition, which can be regarded as giving a ``classical free field realization'' of $\Ocal(\mathfrak{Z}_G)$,
because the chiral analog of $\Ocal(\mathfrak{Z}_G) \rightarrow \Ocal(\mathrm{KT}_G)$ induced by $\varphi$
provides a (chiral) free field realization of the chiral universal centralizer (c.f. section \ref{FFR}).
\begin{prop} [\cite{C, BN}]
\label{immersion}
There exists an open immersion $\varphi : \mathrm{KT}_G \rightarrow \mathfrak{Z}_G.$\\
Moreover, $\varphi$ is symplectic.
\end{prop}

The image of $\varphi$ is described by the open subset $U$, which is the big cell of the Bruhat decomposition of $G$:
$$\operatorname{Im}(\varphi) = (U \times S_{f_0}) \cap  \mathfrak{Z}_G =: \mathfrak{Z}_U.$$

\begin{rmk}
\label{rmk}
The open subset $\mathfrak{Z}_U \subset \mathfrak{Z}_G$ is obtained by using the same method as in the $\mathfrak{Z}_G$ case.
That is, 
\begin{align*}
T^*U \times_{\g^*} S_{f_0}  &\cong U \times S_{f_0}, \\
(U \times S_{f_0})\times_{\g^*} S_{f_0} &\cong \{(g, s) \in U \times S_{f_0} \mid \mathrm{Ad}_g(s)=s \} \\
&= (U \times S_{f_0}) \cap  \mathfrak{Z}_G = \mathfrak{Z}_U.
\end{align*}
\end{rmk}

\subsection{Universal centralizer of $SL_2(\mathbb{C})$}

In this section, we let $G = SL_2(\mathbb{C}).$
Since $$f_0=
\begin{pmatrix}
0 & 0 \\
1 & 0 \\
\end{pmatrix}
, 
h_0=
\begin{pmatrix}
1 & 0 \\
0 & -1 \\
\end{pmatrix}
\text{ and } 
e_0=
\begin{pmatrix}
0 & 1 \\
0 & 0 \\
\end{pmatrix}
,$$ we get
$$S_{f_0} =
\left\{
\begin{pmatrix}
0 & S \\
1 & 0 \\
\end{pmatrix}
\relmiddle|
S \in \mathbb{C}
\right\}.
$$
Therefore, 
\begin{align*}
\mathfrak{Z}_G &=
\left\{
\left(
g
,
\begin{pmatrix}
0 & S \\
1 & 0 \\
\end{pmatrix}
\right)
\in G \times S_{f_0}
\relmiddle|
S \in \mathbb{C}, 
\left[g,
\begin{pmatrix}
0 & S \\	
1 & 0 \\
\end{pmatrix}
\right]
=
0
\right\} \\
&=
\left\{
\begin{pmatrix}
Y & SX \\
X & Y \\
\end{pmatrix}
\relmiddle|
X, Y, S \in \mathbb{C}
, \, \mathrm{det}
\begin{pmatrix}
Y & SX \\
X & Y \\
\end{pmatrix}
=1\right\}
\end{align*}
Hence, $\mathfrak{Z}_G$ is an affine variety in $\mathbb{A}^3_{\mathbb{C}}$ with the coordinate ring
$$\Ocal_{\mathfrak{Z}_G}(\mathfrak{Z}_G) = \frac{\mathbb{C} [X, Y, S]} {( SX^2-Y^2+1 )}.$$
The nontrivial Poisson brackets among the generators of $\Ocal_{\mathfrak{Z}_G}(\mathfrak{Z}_G)$ are as follows:
$$\{S, X\} =Y, \quad \{S, Y \} =SX, \quad \{X,Y \} = -\frac{1}{2} X^2.
$$
Since $U = \{ X \neq 0 \} \subset SL_2(\mathbb{C})$, $\mathfrak{Z}_U = \{ X \neq 0 \} \subset \mathfrak{Z}_G$.
A Poisson algebra $\Ocal_{\mathfrak{Z}_G}(U_X) = \mathbb{C}[X^{\pm 1}, Y]$ is isomorphic to 
$\Ocal_{\mathrm{KT}_G}(\mathrm{KT}_G) = \mathbb{C}[\gamma^{\pm \frac{1}{2}}, b]$ via
$$X \mapsto \gamma^{-\frac{1}{2}}, \quad Y \mapsto -b \gamma^{-\frac{1}{2}}.$$
The restriction map $\, \frac{\mathbb{C} [X, Y, S]} {( SX^2-Y^2+1 )} \rightarrow \mathbb{C}[\gamma^{\pm \frac{1}{2}}, b]$ is given by
\begin{align*}
X &\mapsto \gamma^{-\frac{1}{2}} \, ,\\
Y &\mapsto -b \gamma^{-\frac{1}{2}} \, ,\\
S &\mapsto b^2 + \gamma \, .
\end{align*}
\subsection{The quantized Drinfeld-Sokolov reduction}
Let $\chi : \npm \rightarrow \mathbb{C}$ be a functional on $\npm$ defined by the formula
$$\chi(e^{\alpha} t^n)=
\left\{
\renewcommand{\arraystretch}{1.3}
\begin{array}{ll}
$1$ & \text{(if $\alpha$ is a simple root and $n=-1$)} \\
$0$ & \text{(otherwise).}
\end{array}
\right.
$$
Since $\chi([x,y])=0$ for all $x, y \in \npm$, we get a one dimensional representation of $\npm$, denoted by $\mathbb{C}_{\chi}$.
\begin{defi} [\cite{FF, KRW}]
\label{definition}
Let $M$ be in the category $\Ocal$ of $\gh$ of level $k$.
\\
Define 
$$H_{DS}^{\bullet}(M) := H^{\frac{\infty}{2}+\bullet}(\npm,M \otimes \mathbb{C}_{\chi}),$$
where $H^{\frac{\infty}{2}+\bullet}(\npm, N)$ means the semi-infinite $\npm$-cohomology with coefficients in an $\npm$-module $N$ (\cite{F}).
\end{defi}

We call $H^0_{DS}(M)$ \textbf{the quantized Drinfeld-Sokolov reduction} of $M$. 
It follows from the definition that if $M$ is a vertex operator algebra, then so is $H^0_{DS}(M)$.


\subsection{Chiral differential operators}
\indent
The sheaf of \textbf{chiral differential operators} $\mathcal{D}^{ch}_X$ is the sheaf of conformal vertex algebras on a smooth algebraic variety $X$, defined in \cite{MSV,BD}.
We will concretely describe the structure of the global section of $\mathcal{D}^{ch}_{G}$, where $G$ is an affine algebraic Lie group following \cite{AG}.

We first define \textbf{the arc space} $J_{\infty}X$ of a finite type scheme $X$.
This is a (unique) scheme defined by
$$\operatorname{Hom}_{Scheme}(\operatorname{Spec}(A),J_{\infty}X) \cong \operatorname{Hom}_{Scheme}(\operatorname{Spec}(A \llbracket z \rrbracket),X)$$
for any commutative $\mathbb{C}$-algebra $A$.
If $X =\operatorname{Spec}(R)$ is an affine scheme, $J_{\infty} X = \operatorname{Spec}(J_{\infty} R)$, 
where for $R =\mathbb{C}[x_1, x_2, \ldots , x_n] / (f_1, f_2, \ldots , f_r)$, a differential algebra $(J_{\infty}R, \partial)$ is defined by 
$$J_{\infty}R := \mathbb{C}[\partial^m x_i \mid m \geq 0, 1 \leq i \leq n] / (\partial^l f_j \mid l \geq 0, 1 \leq j \leq r).$$

As a vertex algebra, $\Dch(G) \cong U(\gh_k) \otimes_{U(\g[t] \oplus \mathbb{C}K)} \Ocal(J_{\infty}G),$ where $K$ is the central element in $\gh$ and $J_{\infty}G$ is the arc space of $G$.
Here, $\g[t] \subset J_{\infty}\g$ acts on $\Ocal(J_{\infty}G)$ as a left invariant vector field and $K$ acts as identity.
We have a following vertex algebra structure on the right hand side:\\
we can take $\{x, f \mid x \in \g, f \in \mathcal{O}(G) \}$ as a set of strong generators with fields
$$x(z) = \sum_{n \in \mathbb{Z}} (xt^n) z^{-n-1}, \, f(z) = \sum_{n \in \mathbb{Z}} f_{(n)} z^{-n-1},$$
where $f_{(-n-1)} = \frac{\partial^nf}  {n!}$ if $n \geq 0$ and $0$ if $n<0.$
These are mutually local fields and satisfy the OPEs
$$x(z)y(w) \sim \frac {[x,y](w)} {z-w}+ \frac{k \cdot (x,y)} {(z-w)^2},$$
$$f(z)g(w) \sim 0,$$
$$x(z)f(w) \sim \frac {(x_L f)(w)} {z-w},$$
for $x, y \in \g$ and $f, g \in \Ocal(G)$, where $x_L$ is a left invariant vector field corresponding to $x \in \g.$
Moreover, $U(\gh_k) \otimes_{U(\g[t] \oplus \mathbb{C}K)} \Ocal(J_{\infty}G)$ is a $\mathbb{Z}_{\geq 0}$-graded vertex algebra by setting $\mathrm{deg}(x) = 1, \mathrm{deg}(f) = 0$.

From this fact, the following proposition holds.

\begin{prop}
There are injective vertex algebra homomorphisms 
$$\pi_L : V^k(\g) \rightarrow \Dch(G) \, ; \quad x \mapsto x \otimes 1 \, ,$$
$$j : \Ocal(J_{\infty}(G)) \rightarrow \Dch(G) \, ; \quad f \mapsto 1 \otimes f \, .$$
\end{prop}

In addition, we have a vertex algebra homomorphism $$\pi_R : V^{k^*}(\g) \rightarrow \Dch(G)$$ for $k^* := -k-2h^{\vee}$.
In order to describe this map, let us recall some facts about derivations and differential forms on $G$.

Let $\Omega^1(G)$ be the space of global differential forms on $G$ and $\DerOG$ be the space of derivations on $\Ocal(G).$
Then  we have $\Ocal(G)$-module isomorphisms
$$\Omega^1(G) \xrightarrow{\sim} \mathrm{Hom}_{\mathbb{C}}(\g, \Ocal(G)) \, ; \quad df \mapsto (x \mapsto x_Lf)$$
$$\Ocal(G) \otimes \g \xrightarrow{\sim} \DerOG \, ; \quad 1 \otimes x \mapsto x_L,$$
where $d$ denotes the de Rham differential.

We define $\Omega$ to be the subspace of $\Dch(G)$ spanned by $f \partial g$ with $f,g \in \Ocal(G))$. 
Because of the following lemma, we can identify $\Omega$ with $\Omega^1(G)$.

\begin{lem}
There exists a linear isomorphism
$$\Gamma : \Omega \xrightarrow{\sim} \mathrm{Hom}_{\mathbb{C}}(\g, \Ocal(G)) \, ; \quad f\partial g \mapsto (x \mapsto f \cdot (x_Lg)).$$
\end{lem}

Let $\{ x^i \}_{i = 1}^n$ be a basis of $\g$ and $\{ \omega^i \}_{i = 1}^n$ be the dual basis with respect to the $\Ocal(G)$-bilinear pairing
$\DerOG \times \Omega^1(G) \rightarrow \Ocal(G).$
Since $\{ x^i_L \}_{i = 1}^n$ forms a basis of a free $\Ocal(G)$-module $\Der(\Ocal(G))$,
there exists some invertible matrix $(f^{i,j})$ over $\Ocal(G)$ such that $x^i_R = \sum_{j=1}^n f^{i,j} x^j_L$ for all $i$.

Now we can define $\pi_R$.

\begin{prop} [\cite{AG,GMS}]
\label{commute}
$\mathrm{(i)}$
For $x^i = x^i_{(-1)} \vac \in V^{k^*}(\g),$ set
$$\pi_R(x^i) = \sum_{i =1}^n f^{i, j} x + \sum_{j, k =1}^n k^* \cdot (x^j, x^k) f^{i, j} \omega^k.$$

Then $\pi_R$ gives rise to a vertex algebra embedding $\pi_R : V^{k^*}(\g) \rightarrow \Dch(G)$.\\
$\mathrm{(ii)}$
$\mathrm{Im}(\pi_R) \subset \mathrm{Com}(\mathrm{Im}(\pi_L) , \Dch(G))),$ where $\mathrm{Com}(W, V)$ is the commutant of $W$ in $V$ for a vertex algebra $V$ and its vertex subalgebra $W.$ 
\end{prop}

\begin{rmk}
\label{remark}
We have decreasing (resp. increasing) filtration $F^s$ (resp. $F'_s$) of $\Dch(G)$.
Let $\{a_i\}_{i \in I}$ be a set of strong generators od $\Dch(G).$ Then $F^s$ and $F'_s$ is the space spanned by the vectors
\begin{align*}
&(a_1)_{(-n_1-1)} \cdots (a_m)_{(-n_m-1)} \vac, \, n_1+ \cdots +n_m \geq s, \\
&(a_1)_{(-n_1-1)} \cdots (a_m)_{(-n_m-1)} \vac,  \, \mathrm{deg}(a_1)+ \cdots + \mathrm{deg}(a_m) \leq s,
\end{align*}
respectively.

We have $\gr^F \Dch(G) \cong \gr_{F'} \Dch(G)$ (\cite{A2}).
Let us denote
$$\gr \Dch(G) := \gr^F \Dch(G) \cong \Ocal(\J T^*G) \cong \Ocal(\J \g^*) \otimes \Ocal(\J G).$$ The $\g[t^{\pm 1}]$-action on $\gr \Dch(G)$ is as follows:
\begin{align*}
&xt^n \cdot (a \otimes f_{(m)}) = a \otimes (x_Lf)_{(n+m)}, \\
&xt^{-n-1} \cdot (a \otimes f_{(m)}) = (xt^{-n-1}) a \otimes f_{(m)},
\end{align*}
where $x \in \g$, $n, m \geq 0$, $a \in \Ocal(\J \g^*) \cong \mathrm{Sym}(t^{-1}\g[t^{-1}])$ and $f \in \Ocal(G).$
\end{rmk}


\subsection{Chiral universal centralizer}

The chiral universal centralizer is obtained by performing the quantized Drinfeld-Sokolov reduction of $\Dch(G)$ twice.

There are two commuting actions of $\gh_k$ and $\gh_{k^*}$ on $\Dch(G)$, corresponding to the vertex algebra homomorphisms
\begin{align*}
&\pi_L : V^k(\mathfrak{g}) \rightarrow \Dch(G), \\
&\pi_R : V^{k^*}(\mathfrak{g}) \rightarrow \Dch(G),
\end{align*}
defined in the previous section.

Since $\Dch(G)$ satisfies the condition in Definition \ref{definition}, we can perform the quantized Drinfeld-Sokolov reduction of $\Dch(G)$ first with respect to the former action, and then the latter action.

The first vertex algebra we obtain is said to be \textbf{the equivariant affine $W$-algebra} and denoted by $\mathbf{W}_{G,k}$ (\cite{A1}). That is, $\mathbf{W}_{G,k} := H^0_{DS}(\Dch(G))$.
By Proposition \ref{commute} (ii), we know that $\pi_R$ descends to a vertex algebra homomorphism  $\pi'_R : V^{k^*}(\g) \rightarrow \mathbf{W}_{G,k}$, which gives rise to a $\gh_{k^*}$ action on $\mathbf{W}_{G,k}$.
Then we can consider the quantized Drinfeld-Sokolov reduction of $H^0_{DS}(\Dch(G))$ with respect to this action.
Finally, we get a vertex algebra $$\mathbf{I}_{G,k} := H^0_{DS}(\mathbf{W}_{G,k}),$$ called \textbf{the chiral universal centralizer} associated with $G$.
The level $k$ in $\mathbf{W}_{G,k}$ and $\mathbf{I}_{G,k}$ may be omitted if it is clear.

The chiral universal centralizer is a strict chiral quantization of the (classical)  universal centralizer (\cite{A1}).
That is, $X_{\I} \cong \mathfrak{Z}_G$ as Poisson varieties, $\tilde{X}_{\I} \cong \mathfrak{Z}_G$ as schemes and 
the surjective Poisson vertex algebra homomorphism $\Ocal(\J\tilde{X}_{\I}) \rightarrow \gr \I$ induced by the identity map on $R_{\I}$ is an isomorphism.
Here, for a vertex algebra $V$, $R_V$ denotes the Zhu's $C_2$ algebra of $V$ (\cite{ZhuY}), 
$X_V := \mathrm{Specm}(R_V)$ is the associated variety of $V$ and $\tilde{X}_V := \mathrm{Spec}(R_V)$ is the associated scheme of $V$ (\cite{A1}).

The vertex algebra $\mathbf{I}_{G,k}$ is conformal with central charge $2 \operatorname{rk} (\g) + \langle \rho, \rho^{\vee} \rangle$, 
where $\rho$ (resp. $\rho^{\vee}$) is the half sum of positive root (resp. positive coroot).
Moreover, $\mathbf{I}_{G,k}$ is simple (\cite{AM}) because its associated variety $\mathfrak{Z}_G$ is smooth (\cite{R}) and reduced. 
Note that we have the following vertex algebra homomorphisms:
\begin{align*}
&\pi_L \otimes \pi_R : V^k(\mathfrak{g}) \otimes V^{k^*}(\mathfrak{g}) \rightarrow \Dch(G), \\
&H^0_{DS}(\pi_L) \otimes \pi_R' : \mathcal{W}^k(\mathfrak{g}) \otimes V^{k^*}(\mathfrak{g}) \rightarrow \mathbf{W}_{G,k}, \\
&H^0_{DS}(\pi_L) \otimes H^0_{DS}(\pi_R') : \mathcal{W}^k(\mathfrak{g}) \otimes \mathcal{W}^{k^*}(\mathfrak{g}) \rightarrow \I.
\end{align*}
These maps are conformal when $k$ is non-critical level.

Let the level $k$ be irrational. We have a decomposition of $\mathcal{D}^{ch}_{G,k}(G)$ as a $V^k(\mathfrak{g}) \otimes V^{k^*}(\mathfrak{g})$-module (\cite{AG, ZhuM}):
$$\mathcal{D}^{ch}_{G,k}(G) \cong \bigoplus_{\lambda \in P_+} V^k_{\lambda} \otimes V^{k^*}_{\lambda^*}.$$
Here, $P_+$ is the set of dominant integral weights of $\g$,
$V^k_{\lambda}$ denotes the Weyl module in level $k$ 
and $\lambda^* = -\mathbf{w} \cdot \lambda$, where $\mathbf{w}$ is the longest Weyl group element of $\g$.
It follows that $\mathbf{I}_{G, k}$ has a decomposition as a $\mathcal{W}^k(\mathfrak{g}) \otimes \mathcal{W}^{k^*}(\mathfrak{g})$-module:
$$\mathbf{I}_{G, k} \cong \bigoplus_{\lambda \in P_+} T^k_{\lambda, 0} \otimes T^{k^*}_{\lambda^*, 0},$$
where $T^k_{\lambda, 0} := H^0_{DS}(V^k_{\lambda})$ is the irreducible $\mathcal{W}^k(\mathfrak{g})$-module studied in \cite{AF}.
For $G = SL_2(\mathbb{C})$, this decomposition appeared in \cite{FS} and they called it the modified regular representation of the Virasoro algebra.

\section{Main Results}
\subsection{Free field realization of $\I$}
\label{FFR}
In this section, we construct a free field realization of $\I$ in a similar way to \cite{FS}.

We denote $U$ the big cell of the Bruhat decomposition of $G$, corresponding to the longest Weyl group element $\mathbf{w}$.
This subset $U \subset G$ is open 
 and we can concretely write 
$U=N \widetilde{\mathbf{w}} TN$. Here $\widetilde{\mathbf{w}}$ is the uplift of $\mathbf{w}$ to the normalizer of $T$ in $G$ 
and lies in the connected component of the identity.

The vertex algebra $\Dch(U)$ is equipped with the commuting actions of $\widehat{\mathfrak{g}}_k$ and $\widehat{\mathfrak{g}}_{k^*}$, induced by composing the restriction map with $\pi_L$ and $\pi_R$:
$$V^k(\mathfrak{g}) \xrightarrow{\pi_L} \Dch(G) \xrightarrow{\mathrm{restriction}} \Dch(U),$$
$$V^{k^*}(\mathfrak{g}) \xrightarrow{\pi_R} \Dch(G) \xrightarrow{\mathrm{restriction}} \Dch(U).$$
As in the previous section, we take the quantized Drinfeld-Sokolov reduction twice; first with respect to the $\gh_k$ action, and then with respect to the $\gh_{k^*}$ action.

By performing the first reduction, we obtain the following vertex algebra homomorphisms:
$$ \mathcal{W}^k(\mathfrak{g}) \otimes V^{k^*}(\mathfrak{g}) \rightarrow \mathbf{W}_{G,k} \rightarrow H^0_{DS}(\Dch(U)).$$

Then we perform the second reduction with respect to the homomorphisms
$$V^{k^*}(\mathfrak{g}) \rightarrow \mathbf{W}_{G,k} \rightarrow H^0_{DS}(\Dch(U)),$$
and get the following vertex algebra homomorphisms:
$$\mathcal{W}^k(\mathfrak{g}) \otimes \mathcal{W}^{k^*}(\mathfrak{g}) \rightarrow \I \rightarrow H^0_{DS}(H^0_{DS}(\Dch(U))).$$

Now we state the following theorem, which is the main result of this paper.
\begin{thm}
\label{theorem}
There exists a vertex algebra embedding 
$$\nu: \mathbf{I}_{G,k} \rightarrow \mathcal{D}^{ch}_T(T).$$
Moreover, the corresponding map between associated varieties coincides with $\varphi$ in Prososition \ref{immersion}.
\end{thm}

To prove this theorem, we need the following lemmas.

\begin{lem}
\label{unipotent}
$\mathrm{(i)}$
Let $G$ be a unipotent algebraic Lie group. \\
The action of $t^{-1}\g[t^{-1}]$ on $\Dch(G)$ is free and the action of $\g[t]$ on $\Dch(G)$ is cofree. 
Here, these actions come from the vertex algebra homomorphism $V^k(\g) \rightarrow \Dch(G).$ \\
$\mathrm{(ii)}$
$H^0_{DS}(\mathcal{D}^{ch}_N(N)) \cong \mathbb{C}.$
\end{lem}
\begin{prf}[Lemma \ref{unipotent}]
(i)
Let $M := \Dch(G).$
It suffices to show that $H_{i}(t^{-1}\g[t^{-1}],M)=0$ and $H^{i}(\g[t],M)=0$ if $i \neq 0.$
Set $H_{\bullet}(N_1) := H_{\bullet}(t^{-1}\g[t^{-1}],N_1)$ and $H^{\bullet}(N_2) := H^{\bullet}(\g[t],N_2)$ for a $t^{-1}\g[t^{-1}]$-module $N_1$ and a $\g[t]$-module $N_2$.

We have the $\J \g^*$-action and $\J \g$-action on $\Ocal(\J T^*G)$ (cf. Remark \ref{remark}).
The former action is free and the latter is cofree. Indeed, the freeness is clear.
Since $G$ is unipotent, $\J G$ is so. Hence $\J \g$-action on $\Ocal(\J T^*G)$ is cofree.

Therefore, $H_i(\gr M) = 0$ and $H^i(\gr M) = 0$ for $i \neq 0$, and $H_0(\gr M) \cong \Ocal(\J \g^*)$ and $H^0(\gr M) \cong \Ocal(\J G)$.
We can get the homology (resp. cohomology) spectral sequence $E^r$ (resp. $E_r$) 
such that $E_{\bullet}^1 \cong H_{\bullet}(\gr M)$ and $E_{\bullet}^{\infty} = \gr H_{\bullet}(M)$ (resp. $E^{\bullet}_1 \cong H^{\bullet}(\gr M)$ and $E^{\bullet}_{\infty} = \gr H^{\bullet}(M)$).

Since $H_{i}(\gr M) =0$ if $i \neq 0$ and $H_0(\gr M) \cong \Ocal(\J \g^*)$, we get $E^1=E^{\infty}.$
Therefore, $\gr H_{\bullet}(M) \cong H_{\bullet}(\gr M)$ and we get $H_i(M) = 0$ if $i \neq 0$.
Similarly, $H^i(M) = 0$ if $i \neq 0$.
\qed \\
(ii)
Let $M := \mathcal{D}^{ch}_N(N)$.
Let us consider $C^{\frac{\infty}{2}+\bullet}(\npm,M \otimes \mathbb{C}_{\chi})$, which is the complex that provides 
$H^{\bullet}_{DS}(M) = H^{\frac{\infty}{2}+\bullet}(\npm,M \otimes \mathbb{C}_{\chi})$, as a double complex with derivations coming from the $\n[t^{\pm 1}]$-actions on $M$ and $\mathbb{C}_{\chi}$.
Then we find there exists a convergent cohomology spectral sequence $E'_r$ such that $E'^{\bullet}_1= H^{\frac{\infty} {2} + \bullet}(\n[t^{\pm 1}], M)$ and $E'^{\bullet}_{\infty} = H^{\bullet}_{DS}(M)$.

By Lemma \ref{unipotent} (i) and \cite{V},
$$H^{\frac{\infty} {2} + i}(\n[t^{\pm 1}], M) \cong
\left\{
\renewcommand{\arraystretch}{1.3}
\begin{array}{ll}
\mathrm{Im}(M^{\n[t]} \rightarrow M \rightarrow M/t^{-1}\n[t^{-1}]M) & \text{(if $i=0$)} \\
0 & \text{(if $i \neq 0$).}
\end{array}
\right. $$
Therefore, $E'_r$ collapses at $r =1$ and $E'_1= E'_{\infty}$, so 
$$H^{\frac{\infty} {2} + \bullet}(\n[t^{\pm 1}], M) \cong H^{\bullet}_{DS}(M).$$

We will show that $\mathrm{Im}(M^{\n[t]} \rightarrow M \rightarrow M/t^{-1}\n[t^{-1}]M) = \mathbb{C} \vac.$
It suffices to verify that $M^{\n[t]} = \mathbb{C} \vac \oplus (t^{-1}\n[t^{-1}]M \cap M^{\n[t]})$.

Clearly, $M^{\n[t]} \supset \mathbb{C} \vac \oplus (t^{-1}\n[t^{-1}]M \cap M^{\n[t]})$.
Assume that $M^{\n[t]} \supsetneq \mathbb{C} \vac \oplus (t^{-1}\n[t^{-1}]M \cap M^{\n[t]})$.
Then there exists a nonzero element $a \in M^{\n[t]}$ such that
$$a = 1 \otimes \sum(f_1)_{(-n_1-1)} \cdots (f_m)_{(-n_m-1)} \vac$$
for non-constant functions $f_1, \ldots , f_m \in \Ocal(N)$ and $n_1, \ldots , n_m \geq 0.$
Since $a$ is nonzero, there is an element $xt^n \in \n[t]$ such that $xt^n \cdot a \neq 0.$
This contradicts the assumption that $a \in M^{\n[t]}$.
Therefore, $M^{\n[t]} = \mathbb{C} \vac \oplus (t^{-1}\n[t^{-1}]M \cap M^{\n[t]}).$

Consequently, $\quad H^0_{DS}(\mathcal{D}^{ch}_N(N)) \cong \mathbb{C}.$
\qed
\end{prf}

\vspace{8pt}

\begin{prf}[Theorem \ref{theorem}]
We have $\Dch(U) \cong \mathcal{D}^{ch}_N(N) \otimes \mathcal{D}^{ch}_T(T) \otimes \mathcal{D}^{ch}_N(N)$ as a vertex operator algebra.
The action of $\n[t^{\pm 1}]$ which comes from $\pi_L$ concerns with only the latter $\mathcal{D}^{ch}_N(N)$.
Therefore, we can show 
$$H^0_{DS}(\Dch(U)) \cong \mathcal{D}^{ch}_N(N) \otimes \mathcal{D}^{ch}_T(T),$$
in the same way as Lemma \ref{unipotent} (ii) since $N$ is unipotent.
Likewise, the action corresponding to $\pi_R$ (and $\pi'_R$) concerns with only the former $\Dch(N)$, so
$$H^0_{DS}(H^0_{DS}(\Dch(U))) \cong \mathcal{D}^{ch}_T(T).$$
By composing with $\I \rightarrow H^0_{DS}(H^0_{DS}(\Dch(U))),$
we obtain the map 
$$\nu: \I \rightarrow \mathcal{D}^{ch}_T(T).$$
The injectivity of $\nu$ follows from the simplicity of $\I$.

To complete the proof, we need to show that $\nu$ descends to $\varphi$.
This follows from the commutativity of the Kostant reduction and the quantized Drinfeld-Sokolov reduction (\cite{A4}).
That is, 
\begin{align*}
X_{\mathbf{W}_{G,k}} &\cong X_{\Dch(G)} \times_{\g^*} S_{f_0} \\
&\cong T^*G \times_{\g^*} S_{f_0} \\
&\cong G \times S_{f_0},
\end{align*}
and
\begin{align*}
X_{\I} &\cong X_{\mathbf{W}_{G,k}} \times_{\g^*} S_{f_0} \\
&\cong (G \times S_{f_0}) \times_{\g^*} S_{f_0} \\
&\cong \mathfrak{Z}_G.
\end{align*}
Similarly, we have $X_{\Dch(U)} \cong \mathfrak{Z}_G$ (c.f. Remark \ref{rmk}).
The restriction map $\Dch(G) \rightarrow \Dch(U)$ descends to the inclusion map $T^*U \rightarrow T^*G$,
so $\nu : \I \rightarrow \mathcal{D}^{ch}_T(T)$ descends to $\varphi : \mathrm{KT}_G \cong \mathfrak{Z}_U \rightarrow \mathfrak{Z}_G$.
\qed
\end{prf}

\subsection{An example : $SL_2(\mathbb{C})$-case}
In this section we describe the free field realization $\nu$ concretely  for $G = SL_2(\mathbb{C})$.
In this case, the decomposition
$U=N\mathbf{w}TN$ is described as
$$\left\{
\begin{pmatrix}
1 & x \\
0 & 1 \\
\end{pmatrix}
\begin{pmatrix}
0 & -1 \\
1 & 0 \\
\end{pmatrix}
\begin{pmatrix}
\xi & 0 \\
0 & \xi^{-1} \\
\end{pmatrix}
\begin{pmatrix}
1 & y \\
0 & 1 \\
\end{pmatrix}
\relmiddle|
x,y \in \mathbb{C}, \xi \in \mathbb{C}^{\times}
\right\}.
$$

Let
$e=$
$
\begin{pmatrix}
0 & 1 \\
0 & 0 \\
\end{pmatrix}
$,
$h=$
$
\begin{pmatrix}
1 & 0 \\
0 & -1 \\
\end{pmatrix}
$,
$f=$
$
\begin{pmatrix}
0 & 0 \\
1 & 0 \\
\end{pmatrix}
$
be a basis of $\g$.
Take elements $a, b, c, d$ in $\Ocal(G)$, whose values on $g \in G$ are the $(1, \, 1)$-entry, the $(1, \, 2)$-entry, the $(2, \, 1)$-entry and the $(2, \, 2)$-entry of $g$, respectively.
Then $$\Ocal(G) = \mathbb{C}[a, b, c, d]/(ad-bc-1).$$

We can find
\begin{align*}
&e_R = -d^2e_L-cdh_L+c^2f_L,\\
&h_R = -2bde_L-(ad+bc)h_L+2acf_L,\\
&f_R = b^2e_L +abh_L-a^2f_L,\\
&\omega_e = d \cdot (db)-b \cdot (dd),\\
&\omega_h = \frac{1}{2} (d \cdot (da)+c \cdot (db)-b \cdot (dc)-a \cdot (dd)),\\
&\omega_f =-c \cdot (da) + a \cdot (dc).
\end{align*}
Therefore,
$\Dch(G) \cong U(\gh_k) \otimes_{U(\g[t] \oplus \mathbb{C}K)} \Ocal(J_{\infty}G)$ is strongly generated by 
$e = e_{(-1)} \otimes 1 \, \vac, h = h_{(-1)} \otimes 1 \, \vac, f= f_{(-1)} \otimes 1 \, \vac$ and $a = 1 \otimes a \, \vac, b = 1 \otimes b \, \vac, c = 1 \otimes c \, \vac, d = 1 \otimes d \, \vac$.
More concretely, let $V$ be the vertex algebra generated by fields $e(z), h(z), f(z), a(z), b(z), c(z), d(z)$ with the following nontrivial OPEs:
\begin{align*}
&e(z)h(w) \sim -\frac {2e(w)} {z-w}, \quad
e(z)f(w) \sim \frac {h(w)} {z-w} + \frac{k} {(z-w)^2}, \\
&h(z)h(w) \sim \frac {2k} {(z-w)^2}, \quad
h(z)f(w) \sim \frac {2f(w)} {z-w}, \\
&e(z)b(w) \sim \frac {a(w)} {z-w}, \quad
e(z)d(w) \sim \frac {c(w)} {z-w}, \\
&h(z)a(w) \sim \frac {a(w)} {z-w}, \quad
h(z)b(w) \sim -\frac {b(w)} {z-w}, \\
&h(z)c(w) \sim \frac {c(w)} {z-w}, \quad
h(z)d(w) \sim -\frac {d(w)} {z-w}, \\
&f(z)a(w) \sim \frac {b(w)} {z-w}, \quad
f(z)c(w) \sim \frac {d(w)} {z-w}. 
\end{align*}
Finally, we can obtain $\Dch(G)$ as the quotient of $V$ by the submodule generated by the field $$:a(z)d(z):-:b(z)c(z):-1.$$

Next, we consider $\Dch(U) \cong \mathcal{D}^{ch}_N(N) \otimes \mathcal{D}^{ch}_T(T) \otimes \mathcal{D}^{ch}_N(N)$.

It is easy to see see that $\mathcal{D}^{ch}_N(N)$ is isomorphic to the $\beta \gamma$ system $\langle \beta, \gamma \rangle$, which is the vertex algebra generated by two fields
$\beta(z) = \sum_{n \in \mathbb{Z}} \beta_{(n)} z^{-n-1}, \gamma(z) = \sum_{n \in \mathbb{Z}} \gamma_{(n)} z^{-n}$
and the only nontrivial OPE 
$$\beta(z) \gamma(w) \sim \frac {1} {z-w} .$$
As for $\mathcal{D}^{ch}_T(T)$, this is generated by $\widetilde{h}(z), \lp(z)$ and $\lm(z)$ with the following nontrivial OPE and relation:
$$\widetilde{h}(z)\widetilde{h}(w) \sim \frac {2k} {(z-w)^2}, \quad
\widetilde{h}(z)\lambda_{\pm}(w) \sim \frac {\pm \lambda_{\pm}(w)} {z-w}, $$
$$:\lp(z) \lm(z): -1=0.$$
For simplicity, we choose $\eta(z) := \widetilde{h}(z) -2 :\lm(z) \partial \lp(z):$ instead of $\widetilde{h}(z)$ for the generating field of $\mathcal{D}^{ch}_T(T)$.
That is, $\mathcal{D}^{ch}_T(T)$ is isomorphic to $\el$, the vertex algebra generated by $\eta(z), \lp(z)$ and $\lm(z)$ with the following nontrivial OPE and relation:
$$\eta(z)\eta(w) \sim \frac {2(k+2)} {(z-w)^2}, \quad
\eta(z)\lambda_{\pm}(w) \sim \frac {\pm \lambda_{\pm}(w)} {z-w}, $$
$$:\lp(z) \lm(z): -1=0.$$

Hence, we can write $\Dch(U) \cong \bgx \otimes \el \otimes \bgy$, where $\bgx, \bgy$ are the $\beta	 \gamma$ systems.

The restriction map $r : \Dch(G) \rightarrow \Dch(U)$ is as follows:
\begin{align*}
a(z) \quad &\mapsto \quad :\lp(z) \gx(z): \, , \\
b(z) \quad &\mapsto \quad :\lp(z) \gx(z) \gy(z):- \lm(z) \, , \\
c(z) \quad &\mapsto \quad \lp(z) \, , \\
d(z) \quad &\mapsto \quad :\lp(z) \gy(z): \, , \\
e(z) \quad &\mapsto \quad \by(z) \, , \\
h(z) \quad &\mapsto \quad \eta(z) - 2:\gy(z) \by(z): \, , \\
f(z) \quad &\mapsto \quad -:\lm^2(z) \bx(z):-:\gy^2(z) \by(z):+:\gy(z) \eta(z): +k \partial \gy(z) \, .
\end{align*}

We will take the quantized Drinfeld-Sokolov reduction with respect to the $\gh_k$-action. Then we get $\mathbf{W}_G := H_{DS}^0(\Dch)$ as the quotient of $\widetilde{\mathbf{W}}_G$.
$\widetilde{\mathbf{W}}_G$ is a vertex algebra with generating fields $\a(z), \b(z), \c(z), \d(z), \f(z)$ and nontrivial OPEs
\begin{align*}
&\a(z)\b(w) \sim \frac{:\a^2(w):} {2(z-w)}, \quad
\a(z)\d(w) \sim \frac{:\a(w)\c(w):} {2(z-w)}, \\
&\b(z)\c(w) \sim -\frac{:\a(w)\c(w):} {2(z-w)}, \quad
\c(z)\d(w) \sim \frac{:\c^2(w):} {2(z-w)}, \\
&\b(z)\b(w) \sim \frac{2k+3}{4} \left\{ \frac{:\a^2(w):} {(z-w)^2} + \frac{:(\partial\a)(w)\a(w):} {z-w} \right\}, \\
&\d(z)\d(w) \sim \frac{2k+3}{4} \left\{ \frac{:\c^2(w):} {(z-w)^2} + \frac{:\partial(\c)(w)\c(w):} {z-w} \right\}, \\
&\b(z)\d(w) \sim \frac{2k+3}{4} \cdot \frac{:\a(w)\c(w):} {(z-w)^2} + \frac{2+(2k+3):(\partial\a)(w)\c(w):} {4(z-w)}, \\	
&\f(z)\a(w) \sim \frac {(2k+1)\a(w)} {4(z-w)^2} + \frac {\b(w)} {z-w}, \quad
\f(z)\c(w) \sim \frac {(2k+1)\c(w)} {4(z-w)^2} + \frac {\d(w)} {z-w}, \\
&\f(z)\b(w) \sim - \frac{(k+2)(2k+1)\a(w)} {2(z-w)^3} - \frac{(2k+7)\b(w)} {4(z-w)^2} - \frac{:\f(w)\a(w):} {z-w}, \\
&\f(z)\d(w) \sim - \frac{(k+2)(2k+1)\c(w)} {2(z-w)^3} - \frac{(2k+7)\d(w)} {4(z-w)^2} - \frac{:\f(w)\c(w):} {z-w}, \\
&\f(z)\f(w) \sim -\frac{(k+2)(2k+1)(3k+4)} {2(z-w)^4} - \frac{2(k+2)\f(w)} {(z-w)^2} -\frac{(k+2)(\partial\f)(w)} {z-w}.
\end{align*}
The equivariant affine $W$-algebra $\mathbf{W}_G$ is the quotient of $\widetilde{\mathbf{W}}_G$ by the submodule generated by the field $$:\a(z)\d(z):-:\b(z)\c(z):-:(\partial\a)(z)\c(z):-1.$$

The homomorphism $H^0_{DS}(r) : \mathbf{W}_G \rightarrow \bgx \otimes \el$ is
\begin{align*}
\a(z) \quad &\mapsto \quad :\lp(z) \gx(z): \, , \\
\b(z) \quad &\mapsto \quad -\frac{1}{2}:\gx(z) \eta(z) \lp(z):- \lm(z) \,, \\
\c(z) \quad &\mapsto \quad \lp(z) \,, \\
\d(z) \quad &\mapsto \quad -\frac{1}{2}:\eta(z) \lp(z): \, , \\
\f(z) \quad &\mapsto \quad -\frac{1}{4} :\eta^2(z):-:\bx(z) \lm^2(z):- \frac{k+1}{2} (\partial \eta)(z) \, .
\end{align*}
We will take the quantized Drinfeld-Sokolov reduction with respect to the $\gh_k$-action coming from the vertex algebra homomorphism $V^{k^*}(\g) \rightarrow \Dch(G) \rightarrow \mathbf{W}_G.$

We get $\I$ as the quotient of $\widetilde{\mathbf{I}}_G$, where $\widetilde{\mathbf{I}}_G$ is the vertex algebra generated by $C(z), D(z), F(z)$ subjected to the following nontrivial OPEs:
\begin{align*}
&C(z)D(w) \sim \frac{:C^2(w):} {2(z-w)}, \\
&D(z)D(w) \sim \frac{2k+3}{4} \left\{ \frac{:C^2(w):} {(z-w)^2} + \frac{:(\partial C)(w)C(w):} {z-w} \right\}, \\
&F(z)C(w) \sim \frac{(2k+1)C(w)} {4(z-w)^2} + \frac{D(w)}{z-w}, \\
&F(z)D(w) \sim -\frac{(k+2)(2k+1)C(w)}{(z-w)^3} -\frac{(2k+7)D(w)}{4(z-w)^2}-\frac{:F(w)C(w):}{z-w}, \\
&F(z)F(w) \sim -\frac{(k+2)(2k+1)(3k+4)} {2(z-w)^4} - \frac{2(k+2)F(w)} {(z-w)^2} -\frac{(k+2)(\partial F)(w)} {z-w}.
\end{align*}
The chiral universal centralizer $\I$ is the quotient of $\widetilde{\mathbf{I}}_G$ by the submodule generated by the field 
\begin{align*}
:F(z)C^2(z):+:D^2(z):&-\frac {2k+7} {2}(:C(z) (\partial D)(z): - :(\partial C)(z) D(z):) \\
+ &\frac{2k+7} {4} :(\partial C)^2(z):- \frac{2k+3} {8} :(\partial^2 C)(z) C(z): +1.
\end{align*}

The homomorphism $\nu = H^0_{DS}(H^0_{DS}(r)) : \I \rightarrow \el$ is as follows:
\begin{align*}
C(z) \quad &\mapsto \quad \lp(z), \\
D(z) \quad &\mapsto \quad -\frac{1}{2}:\eta(z) \lp(z):, \\
F(z) \quad &\mapsto \quad -\frac{1}{4} :\eta^2(z):-:\lm^2(z):- \frac{k+1}{2} (\partial \eta)(z).
\end{align*}

When $k =-2$, this coincides with the free field reaization in \cite{BN} by 
$$C \longleftrightarrow X, \quad D \longleftrightarrow -Y, \quad F \longleftrightarrow -S,$$
$$\lambda_{\pm} \longleftrightarrow \gamma^{\mp\frac{1}{2}}, \quad \eta \longleftrightarrow -2b.$$

The conformal vector 
\begin{align*}
:\eta(z) (\partial \lp)(z) \lm(z):&+:\bx(z)(\partial\gx)(z):+:\by(z)(\partial\gy)(z):\\
&+(k+1):(\partial\lp)(z)(\partial\lm)(z):-:(\partial^2\lp)(z)\lm(z):
\end{align*}
in $\Dch(U)$
descends to a conformal vector
\begin{align*}
:\eta(z) (\partial \lp)(z) \lm(z):&-:(\partial\lp)(z)(\partial\lm)(z):\\
&-(k+3):(\partial^2\lp)(z)\lm(z):+(\partial \eta)(z)
\end{align*}
in $\mathbf{I}_{G,k}$ with central charge $26$, by performing the quantized Drinfeld-Sokolov reduction twice.

\renewcommand{\refname}{References}

\end{document}